\documentclass{svmult}

\usepackage{makeidx}         
\usepackage{graphicx}        
\usepackage{multicol}        
\usepackage[bottom]{footmisc}
\usepackage{natbib}
\usepackage{url}

\renewcommand{\PackageWarningNoLine}[2]{}
\newtheorem{algorithm}[theorem]{Algorithm}
\usepackage{amsmath}
\usepackage{amsfonts}
\begin{document}

\title*{{On {A}daptive-{M}ultilevel {BDDC}}}
\author{Bed\v{r}ich Soused\'{\i}k\inst{1,2} \thanks{%
Partially supported by National Science Foundation under grant DMS-0713876
and by the Grant Agency of the Czech Republic under grant 106/08/0403.} \and
Jan Mandel\inst{1} \thanks{%
Supported by National Science Foundation under grant DMS-0713876.}}

\institute{Department of Mathematical and Statistical Sciences, \\
University of Colorado Denver, Campus Box 170, Denver, CO 80217, USA
\and Institute of Thermomechanics, Academy of Sciences of the Czech Republic,
Dolej\v{s}kova 1402/5, 182~00 Prague~8, Czech Republic\\
\texttt{jan.mandel@ucdenver.edu, bedrich.sousedik@ucdenver.edu}
}

\maketitle

\section{Introduction}

\label{mandel_plenary_sec:introduction}

The BDDC method (\cite%
{mandel_plenary_Dohrmann-2003-PSC}) is one of the most advanced methods of iterative substructuring.
In the case of many substructures, solving the coarse problem exactly becomes a bottleneck.
Since the coarse problem has the same structure as the original problem,
it is straightforward to apply the method recursively to solve it only approximately.
The two-level BDDC analysis has been extended into three-levels in a pioneering work
by~\cite{mandel_plenary_Tu-2007-TBT3D, mandel_plenary_Tu-2007-TBT},
and into a general multilevel method by~\cite{mandel_plenary_Mandel-2008-MMB}.
The methods for the adaptive selection of constraints
for the two-level BDDC method have been studied in \cite{mandel_plenary_Mandel-2007-ASF, mandel_plenary_Mandel-2009-ABT}.
Here we combine the two approaches into a new method preserving parallel scalability with increasing
number of subdomains and excellent convergence properties.

The theoretical aspects of the design of the BDDC and a closely related FETI-DP 
on irregular subdomains in the plane has been studied by \cite{mandel_plenary_Klawonn-2008-AFA}.
The authors in particular demonstrated that a proper choice of a certain scaling can significantly improve convergence of the methods.
Our goal here is different. 
We consider only the standard stiffness scaling and we look for a space 
where the action of the BDDC preconditioner is defined. 
A combination of these two approaches, also with the proper choice of initial constraints (\cite{mandel_plenary_Burda-2009-SCB}),
would be of independent interest. 
The presented algorithm has been recently extended into 3D in \cite{mandel_plenary_Sousedik-2010-AMB-thesis}.

\cite{mandel_plenary_Klawonn-2007-IFM,mandel_plenary_Klawonn-2009-HA3} have recently successfully
developed and extensively used several inexact solvers for the FETI-DP method,
and \cite{mandel_plenary_Tu-2008-TBA} has extended the three-level BDDC for the saddle point problems.

All abstract spaces in this paper are finite dimensional. The dual space of
a linear space $U$ is denoted by $U^{\prime }$, and $\left\langle \cdot
,\cdot \right\rangle $ is the duality pairing.

\section{Abstract BDDC for a Model Problem}

\label{mandel_plenary_sec:abstract}Let $\Omega \subset \mathbb{R}^{2}$ be a bounded
polygonal domain, decomposed into $N$ nonoverlapping polygonal substructures
$\Omega _{i}$, $i=1,...,N$, which form a conforming triangulation. That is,
if two substructures have a nonempty intersection, then the intersection is
a vertex, or a whole edge. Substructure vertices will also be called
corners. Let $W_{i}$ be the space of Lagrangean $P1$ or $Q1$ finite element
functions with characteristic mesh size$~h$ on $\Omega _{i}$, and which are
zero on the boundary $\partial \Omega $. Suppose that the nodes of the
finite elements coincide on edges common to two substructures. Let%
\begin{equation*}
W=W_{1}\times \cdots \times W_{N},
\end{equation*}%
and let $U\subset W$ be the subspace of functions that are continuous across
the substructure interfaces. We wish to solve the abstract linear problem
\begin{equation}
u\in U:a(u,v)=\left\langle f,v\right\rangle ,\quad \forall v\in U,
\label{mandel_plenary_eq:problem}
\end{equation}%
for a given $f\in U^{\prime }$, where $a$ is a symmetric positive
semidefinite bilinear form on some space $W\supset U$ and positive definite
on $U$. The form $a\left( \cdot ,\cdot \right) $ is called the energy inner
product, the value of the quadratic form $a\left( u,u\right) $ is called the
energy of $u$, and the norm $\left\Vert u\right\Vert _{a}=a\left( u,u\right)
^{1/2}$ is called the energy norm. The operator $A:U\mapsto U^{\prime }$
associated with $a$ is defined by
\begin{equation*}
a(u,v)=\left\langle Au,v\right\rangle ,\quad \forall u,v\in U.
\end{equation*}

The values of functions from $W$ at the corners and certain averages over the
edges will be called the \emph{coarse degrees of freedom}. Let $\widetilde{W}%
\subset W$ be the space of all functions such that the values of any coarse
degrees of freedom have a common value for all relevant substructures and
vanish on $\partial \Omega $. Define $U_{I}\subset U\subset W$ as the
subspace of all functions that are zero on all substructure boundaries$%
~\partial \Omega _{i}$, $\widetilde{W}_{\Delta }\subset W$ as the subspace
of all function such that their coarse degrees of freedom vanish, $%
\widetilde{W}_{\Pi }$ as the subspace of all functions such that their
coarse degrees of freedom between adjacent substructures coincide and such
that their energy is minimal. Then
\begin{equation}
\widetilde{W}=\widetilde{W}_{\Delta }\oplus \widetilde{W}_{\Pi },\quad
\widetilde{W}_{\Delta }\perp _{a}\widetilde{W}_{\Pi }.  \label{mandel_plenary_eq:tilde-dec}
\end{equation}

The component of the BDDC preconditioner computed in $\widetilde{W}_{\Pi }$
is called the coarse problem, cf.~\cite[Algorithm 11]{mandel_plenary_Mandel-2008-MMB}.
Functions that are $a$-orthogonal to $U_{I}$ are called discrete harmonic.
In \cite{mandel_plenary_Mandel-2005-ATP} and \cite{mandel_plenary_Mandel-2007-ASF}, the analysis was done
in spaces of discrete harmonic functions after eliminating $U_{I}$; this is
not the case here, so $\widetilde{W}$ does not consist of discrete harmonic
functions only. Denote by$~P$ the energy orthogonal projection from $U$ to $%
U_{I}$. Then $I-P$ is known as the projection onto the discrete harmonic
functions.\ Finally, let $E$ be a projection from$~\widetilde{W}$ onto $U$
defined by taking some\ weighted average over the substructure interfaces.

Let us briefly describe the construction of the space $\widetilde{W}$ using the coarse
degrees of freedom. Suppose we are given a space $X$ and a linear operator $%
C:W\mapsto X$ and define
\begin{equation}
\widetilde{W}=\left\{ w\in W:C\left( I-E\right) w=0\right\} .
\label{mandel_plenary_eq:W-tilde}
\end{equation}%
The values $Cw$ will be called the local coarse degrees of freedom. To represent
their common values, i.e. the global coarse degrees of freedom, suppose there is
a space $U_{c}$ and linear operators
\begin{equation*}
Q_{P}^{T}:U\rightarrow U_{c}\quad R_{c}:U_{c}\rightarrow X\quad
R:U\rightarrow W,
\end{equation*}%
such that
\begin{equation*}
CR=R_{c}Q_{P}^{T}.
\end{equation*}%
The operator $Q_{P}^{T}$ selects global coarse degrees of freedom in $U_{c}$
as linear combinations of global degrees of freedom; a global coarse degree
of freedom is given by a row of $Q_{P}$. The operator $R$ (resp. $R_{c}$)\
restricts a vector of global (coarse) degrees of freedom into a vector of
local (coarse) degrees of freedom. See~\cite{mandel_plenary_Mandel-2007-ASF} for more
details.

\subsection{Multilevel BDDC}

\label{mandel_plenary_sec:multilevel}The substructuring components (the domains, spaces and operators) from the previous
section will be denoted by an additional subscript $_{1},$ as $\Omega
_{1}^{i},$ $i=1,\ldots N_{1}$, etc., and called level 1. We will call the
coarse problem in $\widetilde{W}_{\Pi 1}$ the level 2 problem. It has the
same finite element structure as the original problem (\ref{mandel_plenary_eq:problem}) on
level 1, so we have $U_{2}=\widetilde{W}_{\Pi 1}$. Level 1 substructures are
level 2 elements, level 1 coarse degrees of freedom are level 2 degrees of
freedom. The shape functions on level 2 are the coarse basis functions in $%
\widetilde{W}_{\Pi 1}$, which are given by the conditions that the value of
exactly one coarse degree of freedom is one and others are zero, and that
they are energy minimal in $\widetilde{W}_{1}$. Note that the resulting
shape functions on level 2 are in general discontinuous between level 2
elements. Level 2 elements are then agglomerated into nonoverlapping level 2
substructures, etc. Level~$\ell$ elements are level $\ell -1$
substructures, and the level $\ell $ substructures are agglomerates of level
$\ell $ elements. Level $\ell $ substructures are denoted by $\Omega _{\ell
}^{i},$ and they are assumed to form a quasiuniform conforming triangulation
with characteristic substructure size $H_{\ell }$. The degrees of freedom of
level $\ell $ elements are given by level $\ell -1$ coarse degrees of
freedom, and shape functions on level $\ell $ are determined by minimization
of energy on each level $\ell -1$ substructure separately, so $U_{\ell }=%
\widetilde{W}_{\Pi ,\ell -1}$. The averaging operators on level $\ell $, $%
E_{\ell }:\widetilde{W}_{\ell }\rightarrow U_{\ell },$ are defined by
averaging of the values of level $\ell $ degrees of freedom between level $%
\ell $ substructures $\Omega _{\ell }^{i}$. The space $U_{I\ell }$ consists
of functions in $U_{\ell }$ that are zero on the boundaries of all level $%
\ell $ substructures, and $P_{\ell }:U_{\ell }\rightarrow U_{I\ell }$ is the
$a-$orthogonal projection in $U_{\ell }$ onto $U_{I\ell }$. For convenience,
let $\Omega _{0}^{i}$ be the original finite elements, $H_{0}=h$.

\begin{algorithm}[Multilevel BDDC, \protect\cite{mandel_plenary_Mandel-2008-MMB},
Algorithm~17]
\label{mandel_plenary_alg:multilevel-bddc}Define the preconditioner $r_{1}\in U_{1}^{\prime
}\longmapsto u_{1}\in U_{1}$ as follows:

\noindent \textbf{for }$\ell =1,\ldots ,L-1$\textbf{,}

\begin{description}
\item Compute interior pre-correction on level $\ell $,%
\begin{equation}
u_{I\ell }\in U_{I\ell }:a\left( u_{I\ell },z_{I\ell }\right) =\left\langle
r_{\ell },z_{I\ell }\right\rangle ,\quad \forall z_{I\ell }\in U_{I\ell }.
\label{mandel_plenary_eq:ML-uIi}
\end{equation}

\item Get updated residual on level $\ell $,%
\begin{equation}
r_{B\ell }\in U_{\ell },\quad \left\langle r_{B\ell },v_{\ell }\right\rangle
=\left\langle r_{\ell },v_{\ell }\right\rangle -a\left( u_{I\ell },v_{\ell
}\right) ,\quad \forall v_{\ell }\in U_{\ell }.  \label{mandel_plenary_eq:ML-rBi}
\end{equation}

\item Find the substructure correction on level $\ell ,$%
\begin{equation}
w_{\Delta \ell }\in W_{\Delta \ell }:a\left( w_{\Delta \ell },z_{\Delta \ell
}\right) =\left\langle r_{B\ell },E_{\ell }z_{\Delta \ell }\right\rangle
,\quad \forall z_{\Delta \ell }\in W_{\Delta \ell }.  \label{mandel_plenary_eq:ML-wDi}
\end{equation}

\item Formulate the coarse problem on level $\ell $,
\begin{equation}
w_{\Pi \ell }\in W_{\Pi \ell }:a\left( w_{\Pi \ell },z_{\Pi \ell }\right)
=\left\langle r_{B\ell },E_{\ell }z_{\Pi \ell }\right\rangle ,\quad \forall
z_{\Pi \ell }\in W_{\Pi \ell },  \label{mandel_plenary_eq:ML-coarse}
\end{equation}

\item If$\ \ell =L-1$, solve the coarse problem directly and set $%
u_{L}=w_{\Pi L-1}$, \newline
otherwise set up the right-hand side for level $\ell +1$,%
\begin{equation}
r_{\ell +1}\in \widetilde{W}_{\Pi \ell }^{\prime },\quad \left\langle
r_{\ell +1},z_{\ell +1}\right\rangle =\left\langle r_{B\ell },E_{\ell
}z_{\ell +1}\right\rangle ,\quad \forall z_{\ell +1}\in \widetilde{W}_{\Pi
\ell }=U_{\ell +1},  \label{mandel_plenary_eq:ML-ri+1}
\end{equation}
\end{description}

\noindent\textbf{end.}

\noindent \textbf{for }$\ell =L-1,\ldots ,1\mathbf{,}$

\begin{description}
\item Average the approximate corrections on substructure interfaces on
level $\ell $,%
\begin{equation}
u_{B\ell }=E_{\ell }\left( w_{\Delta \ell }+u_{\ell +1}\right) .
\label{mandel_plenary_eq:ML-uBi-1}
\end{equation}

\item Compute the interior post-correction on level $\ell $,%
\begin{equation}
v_{I\ell }\in U_{I\ell }:a\left( v_{I\ell },z_{I\ell }\right) =a\left(
u_{B\ell },z_{I\ell }\right) ,\quad \forall z_{I\ell }\in U_{I\ell }.
\label{mandel_plenary_eq:ML-vIi}
\end{equation}

\item Apply the combined corrections,%
\begin{equation}
u_{\ell }=u_{I\ell }+u_{B\ell }-v_{I\ell }.  \label{mandel_plenary_eq:ML-ui}
\end{equation}
\end{description}

\noindent\textbf{end.}
\end{algorithm}

A condition number bound follows, cf.~[\cite{mandel_plenary_Mandel-2008-MMB}, Lemma~20].

\begin{lemma}
\label{mandel_plenary_lem:bddc-ML-estimate}
If for some $\omega _{\ell }\geq 1$, for all $\ell =1,\ldots ,L-1$,
\begin{equation}
\omega _{\ell }=\sup_{w_{\ell }\in (I-P_\ell)\widetilde{W}_{\ell }}J_\ell(w_\ell), \qquad J_\ell(w_\ell)=\frac{\left\Vert
(I-E_{\ell })w_{\ell }\right\Vert _{a}^{2}}{\left\Vert w_{\ell
}\right\Vert _{a}^{2}},
\end{equation}%
then the multilevel BDDC preconditioner satisfies $\kappa \leq \omega =
{\textstyle\prod_{k=1}^{L-1}}
\omega _{\ell }.$
\end{lemma}

\section{Indicator of the Condition Number Bound}

\label{mandel_plenary_sec:indicator}As in~\cite{mandel_plenary_Mandel-2007-ASF}, we propose as an
indicator of the condition number the maximum of the bounds from Lemma~\ref%
{mandel_plenary_lem:bddc-ML-estimate} computed by considering on each level $\ell$ only
one pair of adjacent substructures $s$ and $t$ at a time:
\begin{equation}
\omega \approx \widetilde{\omega} = \Pi _{\ell =1}^{L-1}\max_{st}\omega
_{\ell }^{st},\quad \omega _{\ell }^{st}=\sup_{w_{\ell }^{st}\in (I-P_{\ell
}^{st})\widetilde{W}_{\ell }^{st}}J_{\ell }^{st}\left( w_{\ell }^{st}\right)
,  \label{mandel_plenary_eq:cond-est}
\end{equation}%
where a pair of substructures is called adjacent if they share an edge, and
the quantities with the superscript $^{st}$ are defined using the domain consisting
of the level $\ell $ substructures $s$ and $t$ only.

The quantity $\widetilde{\omega}$ is called an \emph{indicator} of the condition number bound.

Let $S_{\ell }^{st}$\ be the Schur complement operator
associated with the bilinear form $a\left( \cdot .\cdot \right) $ on the
space $\left( I-P_{\ell }^{st}\right) \widetilde{W}_{\ell }^{st}$. The next
theorem is~\cite[Theorem 2]{mandel_plenary_Mandel-2006-ACS} written in a way suitable for
our purposes.

\begin{theorem}
\label{mandel_plenary_thm:local-est}
Let $a>0$, $\Pi _{\ell }^{st}$ be the orthogonal projection onto $(I-P_\ell^{st})\widetilde{%
W}_{\ell }^{st}$, and $I-\overline{\Pi }_{\ell }^{st}$ the
orthogonal projection onto
\begin{equation*}
\operatorname*{null}\left( \Pi _{\ell }^{st}S_{\ell }^{st}\Pi _{\ell
}^{st}+a\left( I-\Pi _{\ell }^{st}\right) \right) .
\end{equation*}%
Then the stationary values $\omega _{\ell ,1}^{st}\geq \omega _{\ell
,2}^{st}\geq \ldots $ and the corresponding stationary vectors $w_{\ell
,k}^{st}$ of the Rayleigh quotient $J_{\ell }^{st}$ in (\ref{mandel_plenary_eq:cond-est})
satisfy
\begin{equation}
X_{\ell }^{st}w_{\ell ,k}^{st}=\omega _{\ell ,k}^{st}Y_{\ell }^{st}w_{\ell
,k}^{st}  \label{mandel_plenary_eq:local-eig}
\end{equation}%
with $Y_{\ell }^{st}$ positive definite, where \textbf{\ }%
\begin{align*}
X_{\ell }^{st}& =\Pi _{\ell }^{st}\left( I-E_{\ell }^{st}\right) ^{T}S_{\ell
}^{st}\left( I-E_{\ell }^{st}\right) \Pi _{\ell }^{st}, \\
Y_{\ell }^{st}& =\left( \overline{\Pi }_{\ell }^{st}\left( \Pi _{\ell
}^{st}S_{\ell }^{st}\Pi _{\ell }^{st}+a\left( I-\Pi _{\ell }^{st}\right)
\right) \overline{\Pi }_{\ell }^{st}+a\left( I-\overline{\Pi }_{\ell
}^{st}\right) \right) .
\end{align*}%
\end{theorem}

The eigenvalue problem (\ref{mandel_plenary_eq:local-eig}) is obtained by projecting the
gradient of the Rayleigh quotient $J_{\ell }^{st}\left( w_{\ell
}^{st}\right) $ onto the complement in $(I-P_{\ell }^{st})\widetilde{W}%
_{\ell }^{st}$ of the subspace, where its denominator is zero%
, in two steps. Both projections $\Pi _{\ell }^{st}$ and $\overline{\Pi }%
_{\ell }^{st}$ are computed by matrix algebra, which is straightforward to
implement numerically. The projection $\Pi _{\ell }^{st}$ projects onto $%
\operatorname*{null}$ $C_{\ell }^{st}\left( I-E_{\ell }^{st}\right) $, and $I-%
\overline{\Pi }_{\ell }^{st}$ projects onto a subspace of $\operatorname*{null}$ $%
S_{\ell }^{st}$, which can be easily constructed computationally if a matrix
$Z_{\ell }^{st}$ is given such that $\operatorname*{null}$ $S_{\ell }^{st}\subset
\operatorname*{range}Z_{\ell }^{st}$. For this purpose, the rigid body modes are
often available directly or they can be computed from the geometry of the
finite element mesh. For levels $\ell >1$, we can use the matrix $Z_{\ell
}^{st}$ with columns consisting of coarse basis functions, because the span
of the coarse basis functions contains the rigid body modes. In this way, we
can reduce (\ref{mandel_plenary_eq:local-eig}) to a symmetric eigenvalue problem, which is
easier and more efficient to solve numerically.

\section{Optimal Coarse Degrees of Freedom}

\label{mandel_plenary_sec:optimal} Writing ${\widetilde{W}}_{\ell }^{st}{=\operatorname*{null}C}%
_{\ell }^{st}\left( I-E_{\ell }^{st}\right) $ suggests how to add coarse
degrees of freedom to decrease the value of indicator $\widetilde{\omega
}$. The following theorem is an analogy to~\cite[Theorem 3]{mandel_plenary_Mandel-2006-ACS}%
. It follows immediately from the standard characterization of eigenvalues
as minima and maxima of the Rayleigh quotient on subspaces spanned by
eigenvectors, applied to (\ref{mandel_plenary_eq:local-eig}).

\begin{theorem}
\label{mandel_plenary_thm:adaptive}
Suppose $n_{\ell }^{st}\geq 0$ and the coarse dof selection matrix ${%
C_{\ell }^{st}}\left( I-E_{\ell }^{st}\right) $ is augmented by the rows
$w_{\ell ,k}^{stT}\left( I-E_{\ell }^{st}\right) ^{T}%
S_{\ell }^{st}\left( I-E_{\ell }^{st}\right)$, where $%
w_{\ell ,k}^{st}$ are the eigenvectors from (\ref{mandel_plenary_eq:local-eig}). Then $%
\omega _{\ell }^{st}=\omega _{\ell ,n_{\ell }^{st}+1}^{st}$, and $\omega
_{\ell }^{st}\geq \omega _{\ell ,n_{\ell }^{st}+1}^{st}$ for any other
augmentation by at most $n_{\ell }^{st}$ columns.

In particular, if\/ $\omega _{\ell ,n_{\ell }^{st}+1}^{st} \leq \tau$ for all pairs of
adjacent substructures $s,t$ and for all levels $\ell =1,\dots ,L-1$, then $%
\widetilde{\omega }\leq \tau ^{L-1}$.
\end{theorem}

Theorem \ref{mandel_plenary_thm:adaptive} allows us to guarantee that the condition number
indicator $\widetilde{\omega }\leq \tau ^{L-1}$ for a given target value $%
\tau $, by adding the smallest possible number of coarse degrees of freedom.

The primal coarse space selection mechanism that corresponds to this
augmentation can be explained as follows. Let us write the augmentation as
\begin{equation*}
c_{\ell ,k}^{st}=\left[
\begin{array}{cc}
c_{\ell ,k}^{s} & c_{\ell ,k}^{t}%
\end{array}%
\right] =w_{\ell ,k}^{stT}\left( I-E_{\ell }^{st}\right) ^{T}S_{\ell
}^{st}\left( I-E_{\ell }^{st}\right) ,
\end{equation*}%
where $c_{\ell ,k}^{s}$ and $c_{\ell ,k}^{t}$ are blocks corresponding to
substructures $s$ and $t$. It should be noted that the matrix $E_{\ell
}^{st} $\ is constructed for a pair of substructures $s$, $t$\ in such a way
that, cf., e.g.~\cite[eq.~(7)]{mandel_plenary_Mandel-2007-BFM},
\begin{equation*}
B_{D,\ell }^{stT}B_{\ell }^{st}=I-E_{\ell }^{st},
\end{equation*}%
where $B_{D,\ell }^{st}$ and $B_{\ell }^{st}$ are matrices known from the
FETI-DP method. In particular, the entries of $B_{\ell }^{st}$ are $+1$ for
substructure $s$ and $-1$ for substructure $t$. This also relates our
algorithm to the one from~\cite{mandel_plenary_Mandel-2006-ACS}, and~\cite{mandel_plenary_Mandel-2007-ASF}%
. Next, let us observe that, due to the application of $I-E_{\ell }^{st}$,
for the two blocks of $c_{\ell ,k}^{st}$ it holds that $c_{\ell
,k}^{s}=-c_{\ell ,k}^{t}$, i.e., for the two substructures the constraint
weights have the same absolute values and opposite sign. Hence it is
sufficient to consider only one of the two blocks, e.g., $c_{\ell ,k}^{s}$.
The augmentation of the global coarse degrees of freedom selection matrix $%
\left[ Q_{P,\ell },q_{k,\ell }\right] $ is constructed by adding a block of $%
k$ columns computed as
\begin{equation*}
q_{k,\ell }=R_{\ell }^{sT}c_{k}^{sT}.
\end{equation*}%
Each column of $q_{k,\ell }$ defines a coarse degree of freedom associated
with the interface of level $\ell $\ substructures $s$ and $t$. Because $%
R_{\ell }^{s}$ is a $0-1$ matrix, it means that columns in $q_{k,\ell }$ are
formed by a scattering of the entries in $c_{k}^{sT}$.

\section{Adaptive-Multilevel BDDC in 2D}

\label{mandel_plenary_sec:adaptive-multilevel}We describe in more detail the implementation
of the algorithm. It consists of two main steps: (i) setup, and (ii) the
loop of preconditioned conjugate gradients with the Adaptive-Multilevel BDDC
as a preconditioner. The setup was outlined in the previous section, and it
can be summarized as follows:

\begin{algorithm}
Adding of coarse degrees of freedom in order to guarantee that the condition
number indicator $\widetilde{\omega }\leq \tau ^{L-1}$, for a given a target
value$~\tau $:


\bigskip

\noindent \textbf{for levels} $\ell=1:L-1 \mathbf{,}$

\begin{description}
\item \hspace{2.5mm} Create substructures with roughly the same numbers of
degrees of freedom, minimizing the number of \textquotedblleft
cuts\textquotedblright (use a graph partitioner, e.g., METIS~4.0 (\cite%
{mandel_plenary_Karypis-1998-MSP}) with weights on both, vertices and edges).

\medskip

\item \hspace{2.5mm} Find a suitable set of initial constraints (corners in 2D),
and set~up the BDDC structures for the adaptive algorithm
.\end{description}

\textbf{for all edges} $\mathcal{E}_{\ell }$ \textbf{on level} $\ell \mathbf{%
,}$

\begin{description}
\item \hspace{2.5mm} Compute the largest local eigenvalues and corresponding
eigenvectors, until the first $m^{st}$ is found such that $\lambda
_{m^{st}}^{st}\leq \tau $, put $k=1,\dots ,m^{st}$.

\medskip

\item \hspace{2.5mm} Compute the constraint weights ${c}_{k}^{st}=\left[ {c}%
_{k}^{s}\;{c}_{k}^{t}\right] $ as
\begin{equation}
{c}_{k}^{st}=w_{k}^{stT}\Pi _{\ell }^{Ist}\left( I-E_{\ell }^{st}\right)
^{T}S_{\ell }^{st}\left( I-E_{\ell }^{st}\right) \Pi _{\ell }^{Ist},
\label{mandel_plenary_eq:generate constraints}
\end{equation}%
where $\Pi _{\ell }^{Ist}$ is a projection constructed using the set of
initial constraints.

\medskip

\item \hspace{2.5mm} Take one block, e.g., ${c}_{k}^{s}$ and keep nonzero
weights for the edge $\mathcal{E}_{\ell }$.

\medskip

\item \hspace{2.5mm} Add to global coarse dofs selection matrix $Q_{P,\ell }$
the $k$ columns $q_{k,\ell }$ as
\begin{equation}
q_{k,\ell }=R_{\ell }^{sT}c_{k}^{sT}.
\end{equation}
\end{description}

\textbf{end.} 

\begin{description}
\item Setup the BDDC structures for level $\ell$ and check size of the coarse problem: \newline
if small enough, call this the level $L$ problem, factor it directly, and exit from the
loop.
\end{description}

\textbf{end.}
\end{algorithm}

We remark that the adaptive algorithm is significantly simpler and easier to
implement compared to our previous algorithm from~\cite{mandel_plenary_Mandel-2006-ACS},
and~\cite{mandel_plenary_Mandel-2007-ASF}. The constraints in~(\ref{mandel_plenary_eq:generate constraints}%
) are generated from the eigenvectors by the same function that evaluates
the left hand side in~(\ref{mandel_plenary_eq:local-eig}). Then they are \textquotedblleft
torn\textquotedblright\ into two blocks, and entries of one of them, that
correspond to a particular edge shared by substructures $s$\ and $t$ on the
level $\ell $, are scattered into additional columns of the matrix $%
Q_{P,\ell }$.

The adaptive algorithm uses matrices and operators that are readily available in
our implementation of the standard BDDC method (unlike in~\cite{mandel_plenary_Mandel-2009-ABT} this time with an explicit coarse space solve)
with one exception: in order to satisfy the local partition of unity
property, cf.~\cite[eq.~(9)]{mandel_plenary_Mandel-2007-BFM},%
\begin{equation*}
E_{\ell }^{st}R_{\ell }^{st}=I,
\end{equation*}%
we need to generate locally the weight matrices $E_{\ell }^{st}.$

The substructures on higher levels are then treated as (coarse) elements
with energy minimal basis functions. However, the number of added
constraints is a-priori unknown. For this reason, the coarse elements must
allow for variable number of nodes per element, and also for variable number
of degrees of freedom per node. It is also essential to generate a
sufficient number of corners as initial constraints, in particular to
prevent rigid body motions between any pair of adjacent substructures.
This topic has been addressed several times cf., e.g., a recent contribution
by~\cite{mandel_plenary_Burda-2009-SCB}.

Finally, we remark that instead of performing interior pre-correction and
post-correction on level $\ell =1$, cf. eqs (\ref{mandel_plenary_eq:ML-uIi})-(\ref%
{mandel_plenary_eq:ML-rBi}) and (\ref{mandel_plenary_eq:ML-vIi})-(\ref{mandel_plenary_eq:ML-ui}), we can benefit
from reducing the problem to interfaces in the pre-processing step.

\section{Numerical Examples and Conclusion}

\label{mandel_plenary_sec:numerical}The adaptive-multilevel BDDC preconditioner was
implemented in Matlab for the 2D linear elasticity problem (with $\lambda =1$%
, and $\mu =2$) on a square domain discretized by finite elements with $%
1\,182\,722$ degrees of freedom. The domain was decomposed into $2\,304$
subdomains on the second level and into $9$~subdomains on the third-level.
Such decomposition leads to the coarsening ratio $H_{\ell }/H_{\ell -1}=16$,
with $\ell =1,2$. In order to test the adaptive selection of constraints,
one single edge has been jagged on both decomposition levels, see Fig.~\ref%
{mandel_plenary_fig:decomposition}. We have computed the eigenvalues and eigenvectors of (%
\ref{mandel_plenary_eq:local-eig}) by setting up the matrices and using standard methods
for the symmetric eigenvalue problem in Matlab, version 7.8.0.347 (R2009a).

\begin{figure}[tbp]
\centering
\includegraphics[width=11cm]{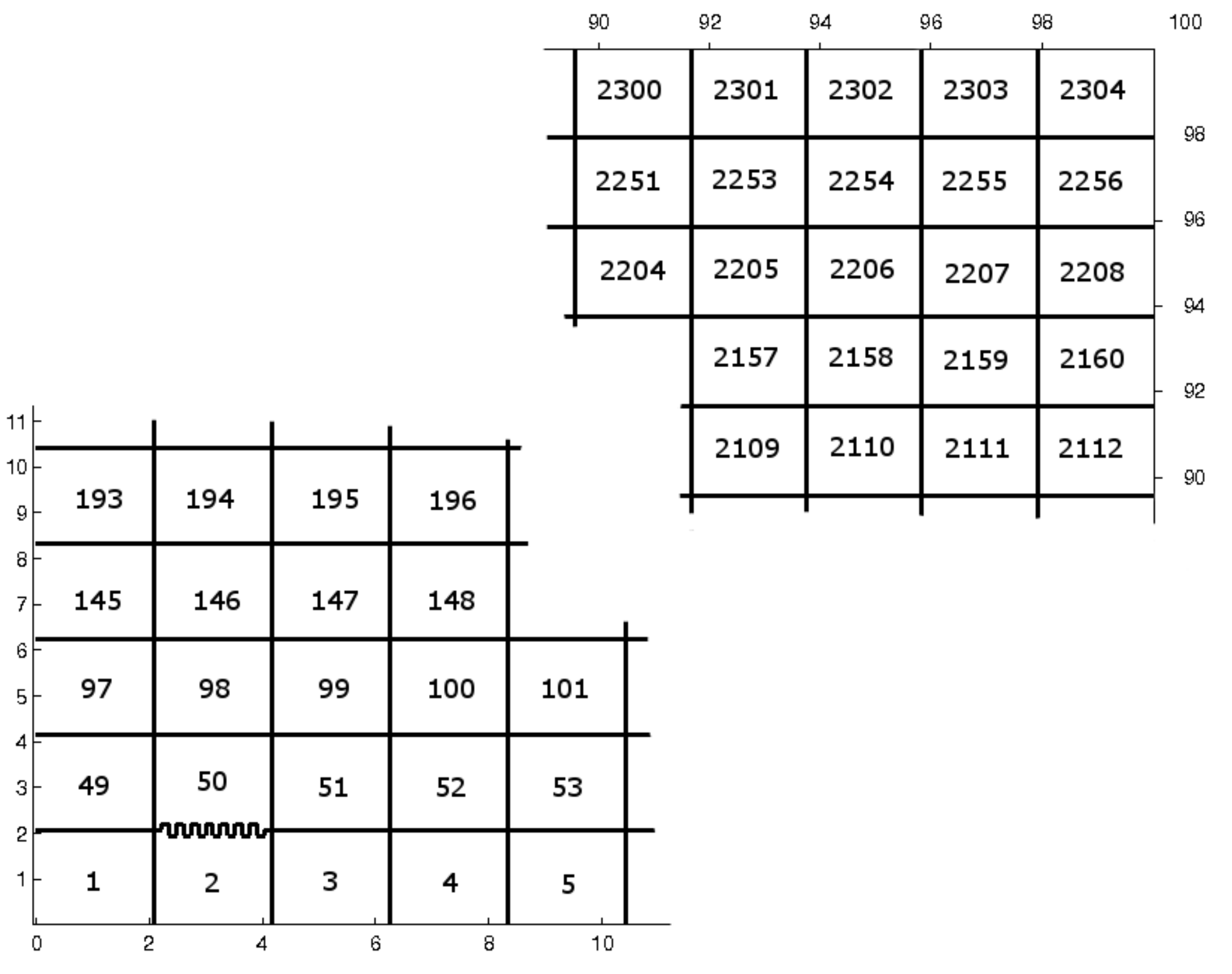}
\includegraphics[width=9cm]{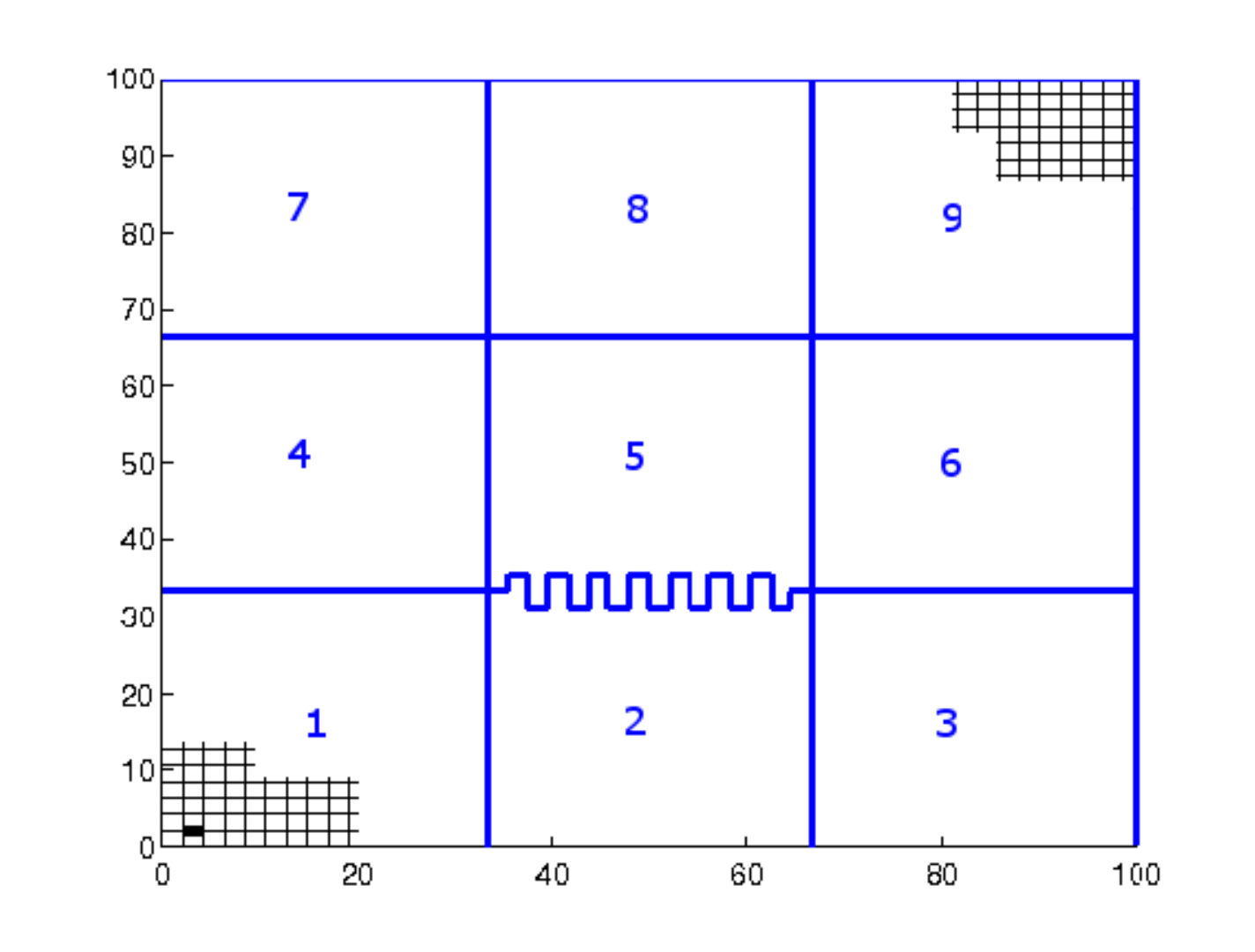}
\caption{The two-level decomposition into $48 \times 48 \, (=2304)$
subdomains (top), and the decomposition into $9$ subdomains for the
three-level method (bottom); the jagged edge from the lower decomposition
level is indicated here by a thick line.}
\label{mandel_plenary_fig:decomposition}
\end{figure}

In the first set of experiments, we have compared performance of the
non-adaptive BDDC method with 2 and 3 decomposition levels. The results are
presented in Tables~\ref{mandel_plenary_tab:2lev-non-adaptive} and~\ref%
{mandel_plenary_tab:3lev-non-adaptive}. As expected from the theory the convergence of the
algorithm deteriorates when additional levels are introduced.

In the next set of experiments, we have tested the adaptive algorithm for
the two-level BDDC. The results are summarized in Table~\ref%
{mandel_plenary_tab:2lev-adaptive}. The algorithm performs consistently with our previous
formulation in~\cite{mandel_plenary_Mandel-2007-ASF}. The eigenvalues associated with edges
between substructures clearly distinguish between the single problematic
edge and the others (Table~\ref{mandel_plenary_tab:2lev-eig}). Adding the coarse dofs
created from the associated eigenvectors according to Theorem~\ref%
{mandel_plenary_thm:adaptive} decreases the value of the condition number indicator $%
\widetilde{\omega }$ and improves convergence at the cost of increasing the
number of coarse dofs.

Finally, we have tested the performance of the Adaptive-Multilevel BDDC for
the model problem with three-level decomposition (Fig.~\ref%
{mandel_plenary_fig:decomposition}). Because the number of coarse degrees of freedom
depends on an a-priori chosen value of$~\tau$ and the coarse basis functions
on level$~\ell$ become shape basis functions on level$~\ell+1$, the
solutions of local eigenvalue problems will depend on$~\tau$ as well.
This fact is illustrated by Table~\ref{mandel_plenary_tab:3lev-eig-tau2} for $\tau=2$,
and by Table~\ref{mandel_plenary_tab:3lev-eig-tau10} for $\tau=10$
(the local eigenvalues for $\tau=3$ were essentially same as for $\tau=2$).
Comparing the values in these two tables, we see that lower values of $\tau$
result in worse conditioning of the local eigenvalue problems on higher decomposition level.
This immediately gives rise to a conjecture that it might not be desirable to decrease the values
of $\tau$ arbitrarily low in order to achieve a better convergence of the method.
On the other hand, for the model problem, comparing the convergence results
for the two-level method (Table~\ref{mandel_plenary_tab:2lev-adaptive})
with the three-level method (Table~\ref{mandel_plenary_tab:3lev-adaptive}), we see that
with the adaptive constraints we were able to achieve nearly the same convergence properties
for the two methods.

\begin{table}[tbp]
\centering
\caption{Results for non-adaptive 2-level method. Constraints are corners,
or corners and arithmetic averages over edges, denoted as c, c+e, resp. $Nc$
is number of constraints, $\mathcal{C}$ is size of the coarse problem
related to size of a subdomain problem, $\protect\kappa$ is the
condition number estimate, $it$ is number of iterations (tol.~$10^{-8}$).}
\label{mandel_plenary_tab:2lev-non-adaptive}
\begin{tabular}{rrrrr}
\hline\noalign{\smallskip}
constraint & $Nc$ & $\mathcal{C}$ & $\kappa$ & $it$ \\
\noalign{\smallskip}\hline\noalign{\smallskip}
c & $4794$ & $9.3$ & $18.41$ & $43$ \\
c+e & $13818$ & $26.9$ & $18.43$ & $32$ \\
\noalign{\smallskip}\hline
\end{tabular}%
\end{table}

\begin{table}[tbp]
\centering
\caption{Results for non-adaptive 3-level method. Headings are as in Table~%
\protect\ref{mandel_plenary_tab:2lev-non-adaptive}.}
\label{mandel_plenary_tab:3lev-non-adaptive}
\begin{tabular}{rrrrr}
\hline\noalign{\smallskip}
constraint & $Nc$ & $\mathcal{C}$ & $\kappa$ & $it$ \\
\noalign{\smallskip}\hline\noalign{\smallskip}
c & $4794 +24$ & $1.0$ & $67.5$ & $74$ \\
c+e & $13818 +48$ & $3.0$ & $97.7$ & $70$ \\
\noalign{\smallskip}\hline
\end{tabular}%
\end{table}

\begin{table}[tbp]
\centering
\caption{Eigenvalues of the local problems for several pairs of subdomains $s
$ and~$t$ on the decomposition level $\ell=1$ (the jagged edge is between
subdomains $2$ and $50$).}
\label{mandel_plenary_tab:2lev-eig}
\begin{tabular}{rrrrrrrrrr}
\hline\noalign{\smallskip}
$s$ & $t$ & $\lambda_{st,1}$ & $\lambda_{st,2}$ & $\lambda_{st,3}$ & $%
\lambda_{st,4}$ & $\lambda_{st,5}$ & $\lambda_{st,6}$ & $\lambda_{st,7}$ & $%
\lambda_{st,8}$ \\
\noalign{\smallskip}\hline\noalign{\smallskip}
$1$ & $2$ & $3.8$ & $2.4$ & $1.4$ & $1.3$ & $1.2$ & $1.1$ & $1.1$ & $1.1$ \\
$1$ & $49$ & $6.0$ & $3.5$ & $2.7$ & $1.4$ & $1.3$ & $1.1$ & $1.1$ & $1.1$
\\
$2$ & $3$ & $5.4$ & $2.6$ & $1.6$ & $1.3$ & $1.2$ & $1.1$ & $1.1$ & $1.1$ \\
$2$ & $50$ & $24.3$ & $18.4$ & $18.3$ & $16.7$ & $16.7$ & $14.7$ & $13.5$ & $%
13.1$ \\
$3$ & $4$ & $3.4$ & $2.4$ & $1.4$ & $1.3$ & $1.1$ & $1.1$ & $1.1$ & $1.1$ \\
$3$ & $51$ & $7.4$ & $4.6$ & $3.7$ & $1.7$ & $1.4$ & $1.3$ & $1.2$ & $1.1$
\\
$49$ & $50$ & $12.6$ & $5.1$ & $4.3$ & $1.9$ & $1.6$ & $1.3$ & $1.2$ & $1.2$
\\
$50$ & $51$ & $8.7$ & $4.8$ & $3.9$ & $1.8$ & $1.5$ & $1.3$ & $1.2$ & $1.2$
\\
$50$ & $98$ & $7.5$ & $4.6$ & $3.7$ & $1.7$ & $1.4$ & $1.3$ & $1.2$ & $1.1$
\\
\noalign{\smallskip}\hline
\end{tabular}%
\end{table}

\begin{table}[tbp]
\centering
\caption{Eigenvalues of the local problems for several pairs of subdomains $s
$, $t$ on level $\ell=2$ with $\protect\tau=2$ (the jagged edge is
between subdomains $2$ and $5$).}
\label{mandel_plenary_tab:3lev-eig-tau2}
\begin{tabular}{rrrrrrrrrr}
\hline\noalign{\smallskip}
$s$ & $t$ & $\lambda_{st,1}$ & $\lambda_{st,2}$ & $\lambda_{st,3}$ & $%
\lambda_{st,4}$ & $\lambda_{st,5}$ & $\lambda_{st,6}$ & $\lambda_{st,7}$ & $%
\lambda_{st,8}$ \\
\noalign{\smallskip}\hline\noalign{\smallskip}
$1$ & $2$ & $16.5$ & $9.0$ & $5.4$ & $2.6$ & $2.1$ & $1.4$ & $1.3$ & $1.3$
\\
$1$ & $4$ & $6.5$ & $4.7$ & $1.9$ & $1.7$ & $1.3$ & $1.2$ & $1.2$ & $1.1$ \\
$2$ & $3$ & $23.1$ & $9.4$ & $4.6$ & $3.2$ & $2.1$ & $1.6$ & $1.4$ & $1.3$
\\
$2$ & $5$ & $84.3$ & $61.4$ & $61.4$ & $55.9$ & $55.8$ & $49.3$ & $48.0$ & $%
46.9$ \\
$3$ & $6$ & $13.7$ & $8.8$ & $4.4$ & $2.2$ & $1.9$ & $1.4$ & $1.3$ & $1.2$
\\
$4$ & $7$ & $6.5$ & $4.7$ & $1.9$ & $1.7$ & $1.3$ & $1.2$ & $1.2$ & $1.1$ \\
$5$ & $6$ & $18.9$ & $13.1$ & $11.3$ & $3.8$ & $2.6$ & $2.1$ & $1.9$ & $1.5$
\\
$5$ & $8$ & $17.3$ & $12.9$ & $10.8$ & $3.6$ & $2.3$ & $2.0$ & $1.8$ & $1.4$
\\
$8$ & $9$ & $13.7$ & $8.8$ & $4.4$ & $2.2$ & $1.9$ & $1.4$ & $1.3$ & $1.2$
\\
\noalign{\smallskip}\hline
\end{tabular}%
\end{table}

\begin{table}[tbp]
\centering
\caption{Eigenvalues of the local problems for several pairs of subdomains $s
$, $t$ on level $\ell=2$ with $\protect\tau=10$ (the jagged edge is
between subdomains $2$ and $5$).}
\label{mandel_plenary_tab:3lev-eig-tau10}
\begin{tabular}{rrrrrrrrrr}
\hline\noalign{\smallskip}
$s$ & $t$ & $\lambda_{st,1}$ & $\lambda_{st,2}$ & $\lambda_{st,3}$ & $%
\lambda_{st,4}$ & $\lambda_{st,5}$ & $\lambda_{st,6}$ & $\lambda_{st,7}$ & $%
\lambda_{st,8}$ \\
\noalign{\smallskip}\hline\noalign{\smallskip}
$1$ & $2$ & $7.7$ & $4.5$ & $2.7$ & $1.6$ & $1.4$ & $1.2$ & $1.2$ & $1.1$ \\
$1$ & $4$ & $3.6$ & $3.0$ & $1.5$ & $1.5$ & $1.2$ & $1.2$ & $1.1$ & $1.1$ \\
$2$ & $3$ & $10.9$ & $4.8$ & $2.7$ & $1.7$ & $1.5$ & $1.2$ & $1.2$ & $1.1$
\\
$2$ & $5$ & $23.2$ & $17.2$ & $13.7$ & $13.7$ & $12.7$ & $12.4$ & $11.0$ & $%
10.9$ \\
$3$ & $6$ & $6.1$ & $4.2$ & $2.5$ & $1.5$ & $1.3$ & $1.2$ & $1.1$ & $1.1$ \\
$4$ & $7$ & $3.6$ & $3.0$ & $1.5$ & $1.5$ & $1.2$ & $1.2$ & $1.1$ & $1.1$ \\
$5$ & $6$ & $9.8$ & $6.2$ & $4.1$ & $2.1$ & $1.6$ & $1.5$ & $1.3$ & $1.2$ \\
$5$ & $8$ & $8.6$ & $5.9$ & $3.9$ & $2.0$ & $1.5$ & $1.4$ & $1.2$ & $1.2$ \\
$8$ & $9$ & $6.1$ & $4.2$ & $2.5$ & $1.5$ & $1.3$ & $1.2$ & $1.1$ & $1.1$ \\
\noalign{\smallskip}\hline
\end{tabular}%
\end{table}

\begin{table}[tbp]
\centering
\caption{Results for the adaptive 2-level method. Headings are same as in Table~%
\protect\ref{mandel_plenary_tab:2lev-non-adaptive}, and $\protect\tau$ is the condition number
target, $\widetilde{\protect\omega}$ is the condition number indicator
.}
\label{mandel_plenary_tab:2lev-adaptive}
\begin{tabular}{rrrrrr}
\hline\noalign{\smallskip}
$\tau$ & $Nc$ & $\mathcal{C}$ & $\widetilde{\omega}$ & $\kappa$ & $it$ \\
\noalign{\smallskip}\hline\noalign{\smallskip}
$\infty$(=c) & $4794$ & $9.3$ & - & $18.41$ & $43$ \\
$10$ & $4805$ & $9.4$ & $8.67$ & $8.34$ & $34$ \\
$3$ & $18110$ & $35.3$ & $2.67$ & $2.44$ & $15$ \\
$2$ & $18305$ & $35.7$ & $1.97$ & $1.97$ & $13$ \\
\noalign{\smallskip}\hline
\end{tabular}%
\end{table}

\begin{table}[tbp]
\centering
\caption{Results for the adaptive 3-level method. Headings are as in Table~%
\protect\ref{mandel_plenary_tab:2lev-adaptive}, but the threshold $\tau$ is now used on each of the two decomposition levels 
and so $\widetilde{\omega} \leq \tau^2$. }
\label{mandel_plenary_tab:3lev-adaptive}
\begin{tabular}{rrrcrr}
\hline\noalign{\smallskip}
$\tau$ & $Nc$ & $\mathcal{C}$ & $\widetilde{\omega}$ & $\kappa$ & $it$ \\
\noalign{\smallskip}\hline\noalign{\smallskip}
$\infty$(=c) & $4794 +24$ & $1.0$ & - & $67.5$ & $74$ \\
$10$ & $4805 +34$ & $1.0$ & $84.97$ & $37.42$ & $60$ \\
$3$ & $18110 +93$ & $3.9$ & $7.88$ & $3.11$ & $19$ \\
$2$ & $18305+117$ & $4.0$ & $3.84$ & $2.28$ & $15$ \\
\noalign{\smallskip}\hline
\end{tabular}%
\end{table}

\bibliographystyle{plainnat}
\bibliography{mandel_plenary}

\end{document}